\documentclass[mathpazo]{cicp}%
\usepackage{amsmath}
\usepackage{amsthm}
\usepackage[T1]{fontenc}
\usepackage[utf8]{inputenc}
\usepackage{units}
\usepackage{subfigure,graphicx}
\usepackage{microtype}
\usepackage{zymacros}
\usepackage{changepage}
\usepackage{url}
\usepackage{booktabs}
\usepackage{amsfonts}
\usepackage{nicefrac}
\usepackage{caption}
\usepackage{xcolor}
\usepackage{subfig}
\usepackage{amssymb}
\usepackage{graphicx}%
\makeatletter

\@ifundefined{showcaptionsetup}{}{ \PassOptionsToPackage{caption=false}{subfig}}
\makeatother

\graphicspath{{./figs/}}
\begin{document}

\title{DeepMartNet - A Martingale based Deep Neural Network Learning Algorithm for
Eigenvalue/BVP Problems and Optimal Stochastic Controls\footnote{Date: 8/23/2023,,
cai@smu.edu. Earlier versions were published in arXiv:2307.11942v2 on 8/3/2023 and  arXiv:2307.11942 on 7/21/2023.
} }
\author{Wei Cai\affil{1}}

\address{
\affilnum{1}\ Dept. of Mathematics, Southern Methodist University, Dallas, TX 75275. \\
}


\begin{abstract}
In this paper, we propose a neural network learning algorithm for solving
eigenvalue problems and boundary value problems (BVPs) for elliptic operators
and initial BVPs (IBVPs) of quasi-linear parabolic equations in high
dimensions as well as optimal stochastic controls. The method is based on the
Martingale property in the stochastic representation for the
eigenvalue/BVP/IBVP problems and martingale principle for optimal stochastic
controls. A loss function based on the Martingale property can be used for
efficient optimization by sampling the stochastic processes associated with
the elliptic operators or value process for stochastic controls. The proposed
algorithm can be used for eigenvalue problems and BVPs and IBVPs with
Dirichlet, Neumann, and Robin boundaries in bounded or unbounded domains and
some feedback stochastic control problems.

\end{abstract}
\date{}
\maketitle

\section{Introduction}

Computing eigenvalue and/or eigenfunctions for elliptic operators or solving
boundary value problem of PDEs, and optimal stochastic control are among the
key tasks for many scientific computing problems, e.g., ground states and band
structure calculation in quantum systems, and financial engineering. Neural
networks have been recently explored for those tasks. FermitNet
\cite{fermitnet} is one of leading methods using anti-symmetrized neural
network wavefunctions in variational Monte Carlo calculation of eigenvalues.
Recently, Han et al \cite{han20} developed a diffusion Monte Carlo method
using the connection between stochastic process and solution of elliptic
equation and the backward Kolmogrov equation to build a loss function for
eigen-value calculations. Based on the same connection, DeepBSDE has also been
designed to solve high dimensional quasi-linear PDEs \cite{DeepBSDE}, which
has also been used for stochastic controls \cite{hanEcontrol}.

In this paper, we use the Martingale problem for the eigenvalue problems and
stochastic controls, and a loss function using fact that the expectation of a
Martingale is constant, thus among any time locations where the expectation
can be approximated by sampling stochastic processes associated with the
elliptic operator and value processes.

\section{DeepMartNet - a Martingale based neural network}

First, we will propose a neural network for computing eigenvalue and
eigenfunction for elliptic operator in high dimensions as arising from quantum
mechanics. It will be apparent that the approach can be applied to solve
boundary value problems of elliptic PDEs and initial boundary value problems
of quasilinear parabolic PDEs.

Consider the following eigenvalue problem%
\begin{align}
\mathcal{L}u+V(\mathbf{x})u  &  =\lambda u,\text{ \ }\mathbf{x}\in D\subset
R^{d},\label{evP}\\
\Gamma(u)  &  =0,\text{ \ }\mathbf{x}\in\partial D,\nonumber
\end{align}
where the boundary operator could be one of the following three cases,%
\[
\Gamma(u)=\left\{
\begin{array}
[c]{cc}%
u & \text{Dirichlet}\\
\frac{\partial u}{\partial n} & \text{Neumann}\\
\frac{\partial u}{\partial n}-cu & \text{Robin}%
\end{array}
\right.  ,
\]
a decay condition will be given at $\infty$ if $D=R^{d},$ and the
differential operator $L$ is given as
\begin{equation}
\mathcal{L}=\mu^{\top}\nabla+\frac{1}{2}Tr(\sigma\sigma^{\top}\nabla \nabla^{\top}%
)\label{gen}%
\end{equation}
and the vector $\mu\mathbf{\in}R^{d},$ matrix $\sigma_{d\times d}$ can be
associated with the drift and diffusion of the following stochastic Ito
process $X_{t}(\omega)$ $\mathbf{\in}$ $R^{d}$ $,\omega\in\Omega$ (random
sample space) with $\mathcal{L}$ as its generator%
\begin{align}
d\mathbf{X}_{t}  &  =\mu dt\mathbf{+}\sigma\mathbf{\cdot}d\mathbf{B}%
_{t}\label{sde}\\
\mathbf{X}_{t}  &  =\mathbf{x}_{0}\in D\nonumber
\end{align}
where $\mathbf{B}_{t}=(B_{t}^{1},\cdots,B_{t}^{d})^{\top}\mathbf{\in}R^{d}$ is
Brownian motion in $R^{d}.$

\subsection{Dirichlet Eigenvalue Problem}

By Ito formula, we have%
\begin{equation}
du(\mathbf{X}_{t})=\mathcal{L}u(\mathbf{X}_{t})dt+\sigma^{\top}\nabla
u(\mathbf{X}_{t})d\mathbf{B}_{t}, \label{ItoP}%
\end{equation}
i.e.,%

\begin{align}
u(\mathbf{X}_{t})  &  =u(x_{0})+\int_{0}^{t}\mathcal{L}u(\mathbf{X}%
_{s})ds+\int_{0}^{t}\sigma^{\top}\nabla u(\mathbf{X}_{s})d\mathbf{B}%
_{s}\nonumber\\
&  =u(x_{0})+\int_{0}^{t}(\lambda-V(\mathbf{X}_{s}))u(\mathbf{X}_{s}%
)ds+\int_{0}^{t}\sigma(\mathbf{X}_{s})^{\top}\nabla u(\mathbf{X}%
_{s})d\mathbf{B}_{s}. \label{ItoForm}%
\end{align}

Due to the fact that the last Ito integral term in (\ref{ItoForm}) is a
Martingale \cite{klebaner} , therefore the following defines Martingale with
respect to a $\mathbf{B}_{t}-$natural filtration $\{\mathcal{F}_{t}\}_{t\geq
0},$%
\begin{equation}
M_{t}=u(\mathbf{X}_{t})-u(\mathbf{x}_{0})-\int_{0}^{t}(\lambda-V(\mathbf{X}%
_{s}))u(\mathbf{X}_{s})ds, \label{Mart}%
\end{equation}
namely, for any $s<t,$%

\begin{equation}
E[M_{t}|\mathcal{F}_{s}]=M_{s}, \label{martin}%
\end{equation}
which implies for any measurable set $A\in$ $\mathcal{F}_{s},$%
\begin{equation}
\int_{A}M_{t}P(d\omega)=\int_{A}E[M_{t}|\mathcal{F}_{s}]P(d\omega
)=\int_{A}M_{s}P(d\omega) \label{Exp1}%
\end{equation}
or%
\begin{equation}
\int_{A}\left(  M_{t}-M_{s}\right)  P(d\omega)=0, \label{Exp2}%
\end{equation}
i.e.,
\begin{equation}
\int_{\Omega}\left(  M_{t}-M_{s}\right)  I_{A}(\omega)P(d\omega)=0
\label{Exp1a}%
\end{equation}
where \ $I_{A}(\omega)$ is the indicator function of the set $A$.

In particular, if we take $A=\Omega\in\mathcal{F}_{s}$ in (\ref{Exp1a}), we have%
\begin{equation}
E[M_{t}-M_{s}]=0. \label{cExp}%
\end{equation}

i.e. the Martingale $M_{t}$ has a constant expectation.

In the case of finite domain $D$, $\tau_{\partial D}$ is a stopping time where
$\tau_{\partial D}$ is the first exit time of the process $X_{t}$ outside $D$,
then $M_{t\wedge\tau_{\partial D}}$ is still a Martingale \cite{klebaner},
thus%
\begin{equation}
E[M_{t\wedge\tau_{\partial D}}-M_{s\wedge\tau_{\partial D}}]=0. \label{cExp1}%
\end{equation}

\begin{remark}
We could define a different generator $\mathcal{L}$ by not including ${\mu}%
^{\top}\nabla$ in (\ref{gen}), then the Martingale in (\ref{Mart1}) will be
changed to
\end{remark}%

\begin{equation}
M_{t}^{\ast}=u(\mathbf{X}_{t})-u(\mathbf{x}_{0})-\int_{0}^{t}(\lambda
-\mu^{\top}(\mathbf{X}_{s})\nabla-V(\mathbf{X}_{s}))u(\mathbf{X}_{s})ds,
\label{Mart1}%
\end{equation}
where the process $\mathbf{X}_{t}$ is given by $d\mathbf{X}_{t}=\sigma
\mathbf{\cdot}d\mathbf{B}_{t},$instead.

\begin{itemize}
\item \bigskip DeepMartNet for eigenvalue $\lambda$
\end{itemize}

Let $u_{\theta}(\mathbf{x})$ be a neural network which will approximate the
eigenfunction with $\theta$ denoting all the weight and bias parameters, for a
given time interval $[0,T]$, we define a partition%
\[
0=t_{0}<t_{1}<\cdots<t_{i}<t_{i+1}<\cdots<t_{N}=T,
\]
and $M$-discrete realizations
\begin{equation}
\Omega^{\prime}=\{\omega_{m}\}_{m=1}^{M}\subset\Omega\label{M-sample}%
\end{equation}
of the Ito process using Euler-Maruyama scheme with $M$-realizations of the
Brownian motions $\mathbf{B}_{i}^{(m)}$, $0\leq m\leq M$, $\mathbf{\ }$%

\[
\mathbf{X}_{i}^{(m)}(\omega_{m})\sim X(t_{i},\omega_{m}),0\leq i\leq N,
\]
where%

\begin{align*}
\mathbf{X}_{i+1}^{(m)}  &  =\mathbf{X}_{i}^{(m)}+\mu\mathbf{(X}_{i}%
^{(m)}\mathbf{)}\Delta t_{i}\mathbf{+}\sigma\mathbf{\mathbf{(X}_{i}%
^{(m)}\mathbf{)}\cdot}\Delta\mathbf{B}_{i}^{(m)},\\
\mathbf{X}_{0}^{(m)}  &  =\mathbf{x}_{0}%
\end{align*}
where $\Delta t_{i}=t_{i+1}-t_{i},$%
\[
\Delta\mathbf{B}_{i}^{(m)}=\mathbf{B}_{i+1}^{(m)}-\mathbf{B}_{i}^{(m)}.
\]

\bigskip

We will build the loss function $l(\theta,\lambda)$ for the eigenfunction
neural network $u_{\theta}(\mathbf{x})$ and the eigenvalue $\lambda$ using the
Martingale property (\ref{Exp2}) and the M-realization of the Ito diffusion
(\ref{sde}).

For each $t_{i}$, we randomly take a subset of $A_{i}\subset\Omega^{\prime}$
with a uniform sampling (without replacement), corresponding to the mini-batch
in computing the stochastic gradient for the empirical training loss, we
should have%
\begin{equation}
\int_{A_{i}}\left(  M_{t_{i+1}}-M_{t_{i}}\right)  P(d\omega)=0,
\end{equation}
which gives an approximate identity for the eigenfunction $u(\mathbf{X}_{t})$
and eigenvalue $\lambda$ using the $A_{i}-$ensemble average,%
\[
\frac{1}{|A_{i}|}\sum_{m=1}^{|A_{i}|}\left(  u(\mathbf{X}_{i+1}^{(m)}%
)-u(\mathbf{X}_{i}^{(m)})-(\lambda-V(\mathbf{X}_{i}^{(m)}))u(\mathbf{X}%
_{i}^{(m)})\Delta t_{i}\right)  \doteq0,
\]
with $|A_{i}|$ being the number of samples in $A_{i}$ (i..e. size of
mini-batch), $\mathbf{X}_{i}^{(m)}=\mathbf{X}_{i}^{(m)}(t_{i},\omega_{m}),$
$\omega_{m}\in A_{i}$, suggesting a loss function, to be used for some epoch(s) of training with 
a given selection of $A_i$'s, in the following form%
\begin{align}
l(\theta,\lambda)  & =l_{\mathbf{x}_{0}}(\theta,\lambda)=\frac{1}{N}\sum
_{i=0}^{N-1}\left(  \frac{1}{|A_{i}|}\sum_{m=1}^{|A_{i}|}\left(  u_{\theta
}(\mathbf{X}_{i+1}^{(m)})-u_{\theta}(\mathbf{X}_{i}^{(m)})-(\lambda
-V(\mathbf{X}_{i}^{(m)}))u_{\theta}(\mathbf{X}_{i}^{(m)})\Delta t_{i}\right)
\right)  ^{2}\label{loss}\\
& +\beta l_{reg}(\theta),\nonumber
\end{align}
where the subscript in $l_{\mathbf{x}_{0}}$ indicates all the sampled paths of
the stochastic process starts from $\mathbf{x}_{0}$, and an regularization
term $l_{reg}(\theta)$ is added for specific needs to be discussed later.

The DeepMartNet approximation for the eigenvalue $\lambda\sim\lambda^{\ast}$
will be obtained by minimizing the loss function $l(\theta,\lambda)$ using
stochastic gradient decent,%
\begin{equation}
(\theta^{\ast},\lambda^{\ast})=\arg\min l(\theta,\lambda). \label{Minimiz}%
\end{equation}

\begin{remark} ({\bf Mini-batch in SGD training and Martingale property})
Due to the equivalence between (\ref{Exp2}) and (\ref{martin}), the loss
function defined above ensures that $M_{t}$ of (\ref{Mart}) for $u_{\theta
}(\mathbf{x})$ will be a Martingale approximately if the mini-batch $A_{i}$
explores all subsets of the sample space $\Omega^{\prime}$ during the SGD optimization process of
the training, and the sample size
$M=|\Omega^{\prime}|$ $\rightarrow$ $\infty,$the time step $\max|\Delta
t_{i}|\rightarrow0,$ and the training converges.

Also, if we take $A_{i}=\Omega^{\prime}$ for all $i,$ there will be no stochasticity
in the gradient calculation for the loss function, we will have a traditional
full gradient descent method and the full Martingale property for  $u_{\theta
}(\mathbf{x})$  is not enforced either. Therefore, the mini-batch practice in DNN SGD optimization
corresponds perfectly with the Martingale definition (\ref{martin}).
\end{remark}

\begin{remark}
(regularizer $l_{reg}(\theta)).$ due to non-uniqueness of the eigenvalues, we
will need to introduce a constrain if we intend to compute the lowest
eigen-value (ground state for quantum systems). The Rayleigh energy can be
used for this purpose for zero drift and constant diffusion coefficient%
\begin{equation}
l_{reg}(\theta)=\int_{\Omega}\left(  \nabla^{\top}u_{\theta}\frac
{{\sigma}{\sigma}^{\top}}{2}\nabla u_{\theta}+Vu_{\theta}%
^{2}\right)  dx+{\gamma}\left(  \int_{\Omega}u_{\theta}^{2}%
d\mathbf{x-}1\right)  ^{2},\label{Regul}%
\end{equation}
where 1-normalization factor for the eigenfunction is also included and the Rayleigh energy integral
can be evaluated with a separate and coarse grid.
\end{remark}

\begin{itemize}
\item \bigskip DeepMartNet for eigenvalue $\lambda$ and eigenfunction $u$
\end{itemize}

As the loss function in (\ref{loss}) only involves paths $\mathbf{X}_{t}$
starting from a fixed point $\mathbf{x}_{0},$ it may not be able to explore
all the state space of the process, therefore the minimization problem in
(\ref{Minimiz}) is expected only to produce a good approximation for the
eigenvalue. To achieve a good approximation to the eigenfunction as well, we
will need to sample the paths of the process $\mathbf{X}_{t}$ from $K$-
initial point $x_{0}^{(k)},1\leq k\leq K,$and define a global loss function%
\begin{equation}
R(\theta,\lambda)=\frac{1}{K}\sum_{k=1}^{K}l_{\mathbf{x}_{0}^{(k)}}%
(\theta,\lambda), \label{loss1}%
\end{equation}
whose minimizer $(\theta^{\ast},\lambda^{\ast})$ is expected to approximate
both the eigenfunction and eigenvalue%
\[
u(x)\sim u_{\theta^{\ast}},\qquad\lambda\sim\lambda^{\ast},
\]
where%

\begin{equation}
(\theta^{\ast},\lambda^{\ast})=\arg\min l(\theta,\lambda). \label{Minimiz1}%
\end{equation}

\subsection{Neumann and Robin Eigenvalue Problem}

We will illustrate the idea for the Robin eigenvalue problem for the simple
case of Laplacian operator,%

\[
\mathcal{L}=\frac{1}{2}\Delta
\]

In probabilistic solutions for Neumann and Robin BVPs, reflecting Brownian
motion will be needed which will go through specular reflections upon hitting
the domain boundary, and a measure of such reflections, the local time of RBM,
will be needed. we will introduce the boundary local time $L(t)$ for
reflecting Brownian motion through a Skorohod problem.

(\textrm{\textbf{Skorohod problem):}} Assume $D$ is a bounded domain in
$R^{d}$ with a $C^{2}$ boundary. Let $f(t)$ be a (continuous) path in $R^{d}$
with $f(0)\in\bar{D}$. A pair $(\xi(t),L(t))$ is a solution to the Skorohod
problem $S(f;D)$ if the following conditions are satisfied:

\begin{enumerate}
\item $\xi$ is a path in $\bar{D}$;

\item $L(t)$ is a nondecreasing function which increases only when $\xi
\in\partial D$, namely,
\begin{equation}
L(t)=\int_{0}^{t}I_{\partial D}(\xi(s))L(ds), \label{sk1}%
\end{equation}

\item The Skorohod equation holds:
\begin{equation}
S(f;D):\qquad\ \xi(t)=f(t)-\int_{0}^{t}n(\xi(s))L(ds), \label{sk2}%
\end{equation}
where $n(x)$ stands for the outward unit normal vector at $x\in\partial D$.
\end{enumerate}

For our case that $f(t)=B_{t}$, the corresponding $\xi_{t}$ will be the
reflecting Brownian motion (RBM) $\mathbf{X}_{t}$. As the name suggests, a RBM
behaves like a BM as long as its path remains inside the domain $D$, but it
will be reflected back inwardly along the normal direction of the boundary
when the path attempts to pass through the boundary. The fact that
$\mathbf{X}_{t}$ is a diffusion process can be proven by using a martingale
formulation and showing that $\mathbf{X}_{t}$ is the solution to the
corresponding martingale problem with the Neumann boundary condition
\cite{hsu84} \cite{Pap1990}.

Due to the fact that RBM $\mathbf{X}_{t}$ is a semimartingale \cite{hsu84}
\cite{Pap1990}, for which the Ito formula \cite{klebaner} give the following%
\begin{equation}
u(\mathbf{X}_{t})=u(x_{0})-\int_{0}^{t}cu(\mathbf{X}_{s})dL(s)-\int_{0}%
^{t}(V(\mathbf{X}_{s})-\lambda)u(\mathbf{X}_{s})ds+\int_{0}^{t}\nabla
u(\mathbf{X}_{s})\cdot d\mathbf{B}_{s}, \label{ItoForm1}%
\end{equation}
where an additional path integral term involving the local time $L(s)$ is
added compared with (\ref{ItoForm}). Again, the last term above being a
Martingale, we can define the following Martingale%
\begin{equation}
M_{t}=u(\mathbf{X}_{t})-u(x_{0})+\int_{0}^{t}cu(\mathbf{X}_{s})dL(s)+\int
_{0}^{t}(V(\mathbf{X}_{s})-\lambda)u(\mathbf{X}_{s})ds. \label{Mart2}%
\end{equation}

Using this Martingale, the DeepMartNet for the Dirichlet eigenvalue problem
can be carried out similarly for the Neumann and Robin eigenvalue problems.
The sampling of reflecting Brownian motion and the computation of local time
$L(t)$ can be found in \cite{ding2023}.

\section{Optimal Stochastic control}

\subsection{Martingale Optimality Principle}

In this section, we will apply the above DeepMartNet
for solving the optimal control of solutions to stochastic differential
equations with a finite time horizon $T$.

Let us consider the following SDE,%
\begin{equation}
d\mathbf{X}_{t}=\mu(t,\mathbf{X}_{t},u_{t})dt\mathbf{+}\sigma(t,\mathbf{X}%
_{t})\mathbf{\cdot}d\mathbf{B}_{t},\text{ \ }0\leq t\leq T \label{ctr-sde}%
\end{equation}
where control $u_{t}$ $\in\mathcal{U}$, where $\mathcal{U}$ is the control
space consisting of $\{\mathcal{F}_{t}\}_{t\geq0}$-predictable processes
taking values in $U\subset R^{m}$ and $\{\mathcal{F}_{t}\}_{t\geq0}$ is the
natural filtration generated by $\mathbf{B}_{t}.$

The running cost of the control problem is a function%
\begin{equation}
c:\Omega\times\lbrack0,T]\times U\rightarrow R, \label{ctr}%
\end{equation}
and for a feedback control, the dependence of $c$ on $\omega\in\Omega$ will be
through the state of the system $\mathbf{X}_{t}(\omega)$, i.e.
\begin{equation}
c(\omega,t,u)=c(\mathbf{X}_{t}(\omega),t,u), \label{fb-ctr}%
\end{equation}
and the terminal cost is defined by a $\mathcal{F}_{T}$-measurable random
variable%
\begin{equation}
\xi(\omega)=\xi(\mathbf{X}_{T}(\omega)) \label{T-cost}%
\end{equation}
where an explicit dependence on the $X_{T}$ is assumed.

For a given control $u$, the total expected cost is then defined by%
\begin{equation}
J(u)=E_{u}[\xi+\int_{[0,T]}c(\mathbf{X}_{t}(\omega),t,u_{t})dt]
\label{total-cost}%
\end{equation}
where the expectation $E_{u}$ is taken with respect to the measure $P^{u}$.

The optimal control problem is to find a control $u^{\ast\text{ }}$such that%
\begin{equation}
u^{\ast\text{ }}=\arg\inf_{u\in\mathcal{U}}J(u). \label{optimalc}%
\end{equation}

To present the Martingale principle for the optimal control, we define the
expected remaining cost for a given control $u$%
\begin{equation}
J(\omega,t,u)=E_{u}[\xi(\mathbf{X}_{T}(\omega))+\int_{[t,T]}c(\mathbf{X}%
_{t}(\omega),t,u_{t})dt|\mathcal{F}_{t}] \label{remain-cost}%
\end{equation}
and a value process%
\begin{equation}
V_{t}(\omega)=\inf_{u\in\mathcal{U}}J(\omega,t,u),\text{ and }E[V_{0}%
]=\inf_{u\in\mathcal{U}} J(u), \label{valueP}%
\end{equation}
and a cost process%
\begin{equation}
M_{t}^{u}(\omega)=\int_{[0,t]}c(\mathbf{X}_{s}(\omega),s,u_{s})ds+V_{t}%
(\omega). \label{costP}%
\end{equation}

The Martingale optimality principle is stated in the following theorem
\cite{cohen15}. \begin{theorem}
(Martingale optimality principle)   $M_{t}^{u}$ is a $P^{u}$-submartingale.
$M_{t}^{u}$ is a $P^{u}$-martingale if and only if control $u=u^{\ast\text{ }%
}$ (the optimal control),and%
\[
E[V_{0}]=E_{u}[M_{0}^{u^{\ast\text{ }}}]=\inf_{u\in\mathcal{U}}J(u).
\]
\end{theorem}

Moreover, the value process $V_{t}(\omega)$ satisfies the following backward
SDE (BSDE)%
\begin{equation}
\left\{
\begin{array}
[c]{c}%
dV_{t}=-H(t,\mathbf{X}_{t},\mathbf{Z}_{t})dt+\mathbf{Z}_{t}dB_{t},0\leq t<T\\
V_{T}(\omega)=\xi(\mathbf{X}_{T}(\omega))
\end{array}
\right.  , \label{Vbsde}%
\end{equation}
where the Hamiltanian%
\[
H(t,\mathbf{x},\mathbf{z})=\inf_{u\in\mathcal{U}}f(t,\mathbf{x,z};u)
\]
and%
\begin{align*}
f(t,\mathbf{x,z};u)  &  =c(\mathbf{x},t,u)+\mathbf{z}\alpha(t,\mathbf{x},u),\\
\alpha(t,\mathbf{x},u)  &  =\sigma^{-1}(t,\mathbf{x})\mu(t,\mathbf{x},u).
\end{align*}

From Pardoux-Peng theory \cite{peng} on the relation between quasi-linear
parabolic equation and backward SDEs , we know that the value process as well
as $Z_{t}(\omega)$ can be expressed in terms of a deterministic function
$v(t,x)$%
\begin{align*}
V_{t}(\omega)  &  =v(t,\mathbf{X}_{t}(\omega))\\
Z_{t}(\omega)  &  =\nabla v(t,\mathbf{X}_{t}(\omega))\sigma(t,\mathbf{X}%
_{t}(\omega))
\end{align*}
where the value function $v(t,\mathbf{x})$ satisfies the following
Hamilton-Jacobi-Bellman (HJB) equation%
\begin{equation}
\left\{
\begin{array}
[c]{c}%
0=\frac{\partial v}{\partial t}(t,\mathbf{x})+\mathcal{L}v(t,\mathbf{x}%
)+H(t,\mathbf{x},\nabla_{x}v\sigma(t,\mathbf{x})),\text{ \ }0\leq
t<T,\mathbf{x}\in R^{d}\\
v(T,\mathbf{x})=\xi(\mathbf{x})
\end{array}
\right.  . \label{hjb}%
\end{equation}

\subsection{DeepMartNet for optimal control $u^{\ast\text{ }}$and value
function $v(t,\mathbf{x})$}

Based on the martingale principle theorem for the optimal feedback control, we
can extend DeepMartNet to approximate the optimal control by a neural network%
\begin{equation}
u_{t}(\omega)=u_{t}(\mathbf{X}(\omega))\sim u_{\theta_{1}}(t,\mathbf{X}%
(\omega)), \label{cdnn}%
\end{equation}
where $u_{\theta_{1}}(t,\mathbf{x})\in C([0,T]\times R^{d}$ ) will be a neural
network approximation for a $d+1$ dimensional function with network parameters
$\theta_{1}$, and the value function by another network%
\begin{equation}
v(t,\mathbf{x})\sim v_{\theta_{2}}(t,\mathbf{x}). \label{vdnn}%
\end{equation}

The loss function will consist of two parts, one for the control network and
one for the value network%
\[
l(\theta_{1},\theta_{2})=l_{ctr}(\theta_{1})+l_{val}(\theta_{2})
\]
where, similar to (\ref{loss}),%
\begin{align}
l_{ctr}(\theta_{1}) &  =l_{ctr,\mathbf{x}_{0}}(\theta_{1})\nonumber\\
&  =\frac{1}{N}\sum_{i=0}^{N-1}\frac{1}{|A_{i}|}\sum_{m=1}^{|A_{i}|}\left(
c(X_{t_{i}},t_{i},u_{\theta_{1}}(t_{i},\mathbf{X}_{i}^{(m)}))\Delta
t_{i}+v_{\theta_{2}}(t_{i+1,}\mathbf{X}_{i+1}^{(m)})-v_{\theta_{2}}%
(t_{i,}\mathbf{X}_{i}^{(m)})\right)  ^{2}\label{loss-ctr}%
\end{align}
and, by using Ito formula for $v_{\theta_{2}}(t,\mathbf{x})$, we can obtain a
similar Martingale form for the HJB equation (\ref{hjb}) and define a similar
loss function for the value function $v(t,\mathbf{x})$ as in (\ref{loss})%
\begin{align}
&  l_{val}(\theta_{2})=l_{val,\mathbf{x}_{0}}(\theta_{2})\nonumber\\
&  =\frac{1}{N}\sum_{i=0}^{N-1}\left(  \frac{1}{|A_{i}|}\sum_{m=1}^{|A_{i}%
|}\left(
\begin{array}
[c]{c}%
v_{\theta_{2}}(t_{i+1,}\mathbf{X}_{i+1}^{(m)})-v_{\theta_{2}}(t_{i,}%
\mathbf{X}_{i}^{(m)})+\\
H(t_{i},\mathbf{X}_{i}^{(m)},\nabla_{x}v_{\theta_{2}}(t_{i,}\mathbf{X}%
_{i}^{(m)})\sigma(t,\mathbf{X}_{i}^{(m)}))\Delta t_{i}%
\end{array}
\right)  \right)  ^{2}\label{loss-val}\\
&  +\beta\frac{1}{M}\sum_{m=1}^{M}(v_{\theta_{2}}(T,\mathbf{X}_{N}^{(m)}%
)-\xi(\mathbf{X}_{N}^{(m)}))^{2}.\nonumber
\end{align}

Again, for better accuracy globally for the control and value networks, we can
define a global loss function with more sampling of the starting points
$\mathbf{x}_{0}^{(k)},1\leq k\leq K$,%
\begin{equation}
R(\theta_{1},\theta_{2})=\frac{1}{K}\sum_{k=1}^{K}\left(  l_{ctr,\mathbf{x}%
_{0}^{(k)}}(\theta_{1})+l_{val,\mathbf{x}_{0}^{(k)}}(\theta_{2})\right)  .
\end{equation}

The above approach requires an accurate result for the value function
$v(t_{i,}\mathbf{X}_{i}^{(m)})$ in the region explored by the process
$\mathbf{X}_{i}^{(m)},$ this could pose a challenge to the DeepMartNet. An
alternative approach is to use FBSDE based learning algorithm in
\cite{zhang22}, which has been shown to be able to meet this requirement.

\section{Conclusion}

In this paper, we introduce a Martingale based neural network for finding the
eigenvalue and eigenfunction for general elliptic operators for general types
of boundary conditions, solutions of BVPs and IBVPs of PDEs as well as optimal
stochastic controls. Future numerical experiments will be carried out to
evaluate the efficiency and accuracy of the proposed algorithm, especially in
high dimensions.

\end{document}